\documentclass{amsart}
\usepackage{amsfonts,latexsym,graphicx,amssymb,amsmath,url,amsthm}
\newtheorem{theorem}{Theorem}
\newtheorem{lemma}{Lemma}

\newtheorem{conjecture*}{Conjecture}
\theoremstyle{definition}
\theoremstyle{remark}
\newtheorem{remark}{{\bf{Remark}}}

\newcommand{\Z}{\mathbb Z}
\newcommand{\Q}{\mathbb Q}

\newcommand{\F}{\mathbb F}

\begin{document}

\title{Systems of cubic forms over $p$-adic fields}
\author[R. Dietmann] {Rainer Dietmann}
\address{Department of Computing, Security and Mathematics, Royal Holloway,
University of London, Egham, Surrey, TW20 0EX, United Kingdom}
\email{rainer.dietmann@rhul.ac.uk} 

\begin{abstract}
We show that if $p \ge 5$ and $C_1, \ldots, C_r \in \Q_p[X_1, \ldots, X_n]$
are cubic forms with $n>72r^2 \frac{1}{1-\log 2/\log p}$, then the
system $C_1=\ldots=C_r=0$ has a common non-trivial $p$-adic zero.
Our approach first follows the author's previous work based on reduction
theory, Warning's second theorem and an upper
bound for the number of singular
zeros of a cubic form over a finite field in terms of its $h$-invariant.
In this key lemma we improve the saving in the
exponent from order of magnitude $h^{1/2}$ to order of magnitude
$h$.
To this end we apply recent work by Kazhdan and Polishchuk, and
Kazhdan, Lampert and Polishchuk regarding the $h$-invariant under
field extensions and its connection to the singular locus, where we
also make use of
B\'{e}zout's theorem in the form of work by Heath-Brown on systems of
quadratic forms. 
\end{abstract}

\maketitle

\section{Introduction}
Let $p$ be a rational prime. A classical result by Lewis and Dem'yanov
\footnote{Dem'yanov assumed $p \ne 3$.}
(see \cite{Lewis}, \cite{Demyanov})
asserts that every cubic form $C(X_1,\ldots,X_n) \in
\Q_p[X_1, \ldots, X_n]$ with $n \ge 10$ has a non-trivial
$p$-adic zero, the latter bound $n \ge 10$ being sharp in general.
More generally, by a conjecture of Artin a system
of $r$ cubic forms in $n$ variables
over $\Q_p$ should have a
common non-trivial $p$-adic zero provided that $n \ge 9r+1$. Indeed,
work of Ax and Kochen \cite[Theorem 5]{AK} shows that $n \ge 9r+1$
suffices if in addition $p$ exceeds a bound depending on $r$, but
a full proof of Artin's conjecture for systems of cubic forms for
$r>1$ doesn't seem to be in the cards at present.
\footnote{Though predictions are always difficult, especially about
the future, and even more so
given recent AI advances in mathematics.}
The best result for large $r$, valid uniformly for all primes $p$,
is due to Schmidt \cite[Theorem 1]{Schmidt},
who obtained a sufficient lower bound $n>5300 r (3r+1)^2$ on the
number of variables. His approach uses cubic exponential sums.
Using a more elementary approach, Leep and Schmidt, and Wooley
(see \cite{LW} and \cite{Wooley}) obtained
bounds that scale like $r^4$ but are superior to Schmidt's result for
small values of $r$. The author \cite[Theorem 1]{RD}
developed an approach based on
reduction theory and quadratic exponential sums to prove a key lemma
(see \cite[Lemma 5]{RD})
bounding the number of singularities of a cubic form over a finite field
in terms of its $h$-invariant,
which works for $p \ne 2$ and leads to a cubic bound of the form
$n \ge 125 r^3 + O(r^2)$.
In this note we want to make use of
recent advances regarding the $h$-invariant of cubic forms
(also known as `slice rank' or `Schmidt rank')
and its connection to the `Birch rank', see \cite{KP}, \cite{KAP}.
For primes $p$ at least $5$ this allows us to replace
the aforementioned key lemma in \cite{RD}
by the more efficient Lemma \ref{key} in this paper which leads to a result
scaling quadratically in the number $r$ of cubic forms.

\begin{theorem}
Let $p$ be a rational prime at least $5$, and let $C_1, \ldots, C_r
\in \Q_p[X_1, \ldots, X_n]$ be cubic forms with
\[
  n>\frac{72r^2}{1-\log 2/\log p}.
\]
Then the system $C_1=\ldots=C_r=0$ has a non-trivial $p$-adic zero.
\end{theorem}

\noindent
Note that since
\[
  \frac{72}{1-\log 2/\log 5}<127,
\]
the lower bound $n \ge 127r^2$ works for every $p \ge 5$.

\bigskip
\noindent
\textbf{AI declaration:} The key ideas in this paper were all found without
the use of any AI, and the preceding main result and its proof were
presented at several conferences (`Diophantine Equations, Combinatorics,
Analysis in Number Theory', ICMS, June 2025; `Full circle: 100 years
of the circle method', IML, July 2025; `ELAZ 2026', February 2026).
The author wrote the present manuscript (after ChatGPT wrote a
first alternative draft from the author's ELAZ 2026 slides), but
ChatGPT 6.1/6.2 Astra Pro provided some assistance in exposition,
proofreading and proof checking.
Moreover, after being prompted with the author's idea of using B\'{e}zout's
theorem in the form of Lemmas 5 and 6 in \cite{HB}, it provided
the streamlined proof of Lemma \ref{bezout}.

\section{The $h$-invariant and its connection to singularities}
For a field $K$ and a cubic form $C \in K[X_1, \ldots, X_n]$, we define
its \emph{$h$-invariant $h_K(C)$} to be the smallest non-negative
integer $h$ such that we can write
\begin{equation}
\label{easy}
  C=L_1 Q_1 + \ldots + L_h Q_h
\end{equation}
for suitable linear and quadratic forms $L_i, Q_i$, respectively,
with coefficients in $K$. It is easy to see that $n-h_K(C)$ is the
maximal dimension of any $K$-linear space on which $C$ vanishes
identically and
this property is often crucial for arithmetic applications.
Writing $\overline{K}$ for the algebraic closure of a field
$K$, we have the trivial inequality $h_{\overline{K}}(C) \le h_K(C)$,
but inequalities in the other direction are harder to establish and
provide a link from the arithmetic/algebraic side to the geometric
world:
Recent work of Kazhdan and Polishchuk (see \cite[Theorem A]{KP} with
$d=3$) shows that if
$K$ is a perfect field and $C \in K[X_1, \ldots, X_n]$ a cubic form, then
\begin{equation}
\label{magic1}
  h_K(C) \le 3 h_{\overline{K}}(C).
\end{equation}
In particular, this holds true for cubic forms $C$ with coefficients
in a finite field $K=\F_p$.
The quantity $h_{\overline{K}}(C)$ can now be linked to the
dimension of the singular locus of the hypersurface $C=0$.
Again, there is an easy observation, namely that by
setting $L_1= \ldots = L_h=Q_1= \ldots =Q_h=0$ in \eqref{easy}
one can produce an
algebraic set of singular zeros of $C$ and this way bound above the
codimension of the singular locus of $C=0$ in terms of $h_K(C)$.
Regarding the harder problem of an inequality in the reverse direction,
recent work of Kazhdan, Lampert and
Polishchuk \cite[Theorem 1.3 (i) with $\mathbf{k}=\overline{\F}_p$
and $d=3$]{KAP} shows that if $p \ge 5$ and
$C  \in \F_p[X_1, \ldots, X_n]$ is a cubic form, then
\begin{equation}
\label{magic2}
  h_{\overline{\F}_p}(C) \le 2 \; \text{codim}
  \{\mathbf{x} \in \overline{\F}_p^n: \nabla C(\mathbf{x})= \mathbf{0}\}.
\end{equation}
For a cubic form $C$ in $n$ variables,
the condition $\nabla(C) = \mathbf{0}$ leads
to a system of $n$ homogeneous quadratic forms simultaneously
being zero.
We need an upper bound  for the number of points on such algebraic sets
over a finite field.
\begin{lemma}
\label{bezout}
Let $p$ be a rational prime,
and let $V \subset \overline{\F}_p^n$ be a non-empty
algebraic set cut out by polynomials of total degree at most two.
If $\text{codim}(V)=s$, then
\[
  \#(V \cap \F_p^n) \le 2^s p^{n-s}.
\]
\end{lemma}
\begin{proof}
Decompose
\[
  V = V_1 \cup \ldots \cup V_t
\]
into distinct geometric irreducible components $V_i$, and write
\[
  d_i = \text{dim} V_i, \quad \delta_i=\text{deg} V_i \quad (1 \le i \le t).
\]
Then by Lemma 5 in \cite{HB} we have
\begin{equation}
\label{hb1}
  \#(V_i \cap \F_p^n) \le \delta_i p^{d_i} \quad (1 \le i \le t).
\end{equation}
Moreover, by an application of B\'{e}zout's theorem,
see Lemma 6 in \cite{HB} (putting $d=2$ there), we obtain
\begin{equation}
\label{hb2}
  \sum_{i=1}^t \delta_i 2^{d_i} \le 2^n.
\end{equation}
As $d_i \le n-s \quad (1 \le i \le t)$
from \eqref{hb1} and \eqref{hb2} we obtain
\[
  \#(V \cap \F_p^n) \le \sum_{i=1}^t \delta_i p^{d_i}
  \le \left( \frac{p}{2} \right)^{n-s} \sum_{i=1}^t \delta_i 2^{d_i}
  \le 2^s p^{n-s}
\]
as required.
\end{proof}
\begin{remark} The previous lemma is true but unfortunately
trivial for $p=2$.
\end{remark}

We are now ready to state and prove our new key lemma.
\begin{lemma}
\label{key}
Let $p \ge 5$ and $C \in \F_p[X_1, \ldots, X_n]$ be a cubic form
with $h_{\F_p}(C) \ge h$, where $h$ is a non-negative real number. Then
\[
  \#\{\mathbf{x} \in \F_p^n: \nabla C(\mathbf{x})=\mathbf{0}\}
  \le 2^{h/6} p^{n-h/6}.
\]
\end{lemma}
\begin{proof}
Let $V \subset \overline{\F}_p^n$ be the algebraic set defined by the
vanishing of all partial derivatives $\frac{\partial C}{\partial X_i} \quad
(1 \le i \le n)$, which are quadratic forms. Then
\[
  \{\mathbf{x} \in \overline{\F}_p^n: \nabla C(\mathbf{x})=\mathbf{0}\}
  =V
\]
and therefore
\[
  \#\{\mathbf{x} \in \F_p^n: \nabla C(\mathbf{x})=\mathbf{0}\}
  = \# (V \cap \F_p^n).
\]
Now as $h_{\F_p}(C) \ge h$, from \eqref{magic1} and \eqref{magic2} we obtain
\[
  s = \text{codim} V \ge \frac{h_{\F_p}(C)}{6} \ge \frac{h}{6},
\]
whence Lemma \ref{bezout} and the observation that $2^s p^{n-s}$
is a decreasing function in $s$ furnish the upper bound
\[
  \# (V \cap \F_p^n) \le 2^s p^{n-s} \le 2^{h/6} p^{n-h/6}
\]
as required.
\end{proof}

Lemma \ref{key} is our new weapon: Our previous version
\cite[Lemma 5]{RD} had a saving of order of magnitude $h^{1/2}$ in
the exponent rather than our linear saving $h/6$ now.

\section{Proof of the main result}
Without loss of generality we can assume that the cubic forms
$C_i$ have coefficients in $\Z_p$, and by the compactness argument
preceding \cite[Lemma 1]{RD}, it suffices to treat a
`reduced' system $(C_1, \ldots, C_r)$.
Write $*$ for reduction modulo $p$. The system
$(C_1, \ldots, C_r)$ being reduced implies that for all
$\mathbf{a} \in \F_p^r \backslash \{\mathbf{0}\}$ we have
\begin{equation}
\label{h}
  h_{\F_p}(a_1 C_1^* + \ldots + a_r C_r^*) \ge \frac{n}{3r}
\end{equation}
(see \cite[Lemma 1]{RD}). Following the framework laid out in \cite{RD}, our
strategy is to find an $\mathbf{x} \in \F_p^n$ for which
\[
  C_1^*(\mathbf{x}) = \ldots = C_r^*(\mathbf{x}) = 0
\]
and
\[
  \nabla C_1^*(\mathbf{x}), \ldots, \nabla C_r^*(\mathbf{x}) 
\]
are linearly independent. Such a non-singular zero $\mathbf{x}$ modulo $p$
can then via Hensel's lemma be lifted to a common $p$-adic zero of
$C_1, \ldots, C_r$. To this end we use a counting argument: By
`Warning's  second theorem' \footnote{which is actually Satz 3 (!) in
\cite{Warning}, together with the remark `Entsprechendes gilt f\"ur
die simultanen L\"osungen mehrerer Polynomgleichungen:', see pages
79 and 82.}
applied to the special case of a system
of $r$ cubic forms over $\F_p$, we have
\begin{equation}
\label{auf_dem}
  \#\{\mathbf{x} \in \F_p^n: C_1^*(\mathbf{x}) = \ldots = C_r^*(\mathbf{x})
  = 0\} \ge p^{n-3r}.
\end{equation}
Having established a lower bound for the number of
all common zeros of $C_1^*, \ldots, C_r^*$, we need an upper bound for 
the number of singular zeros: If $\nabla C_1^*(\mathbf{x}), \ldots,
\nabla C_r^*(\mathbf{x})$ are linearly dependent for
some $\mathbf{x} \in \F_p^n$, then there exists
$\mathbf{a} \in \F_p^r \backslash \{\mathbf{0}\}$ such that
\[
  a_1 \nabla C_1^*(\mathbf{x}) + \ldots +
  a_r \nabla C_r^*(\mathbf{x})=\mathbf{0}.
\]
For
\[
  C_{\mathbf{a}} = a_1 C_1^* + \ldots + a_r C_r^*
\]
by linearity of $\nabla$ we then have
\[
  \nabla C_{\mathbf{a}}(\mathbf{x})=\mathbf{0}.
\]
Now from \eqref{h} we have
\[
  h_{\F_p}(C_{\mathbf{a}}) \ge \frac{n}{3r}.
\]
Hence by Lemma \ref{key} we obtain
\[
  \#\{\mathbf{x} \in \F_p^n: \nabla C_{\mathbf{a}}(\mathbf{x})=\mathbf{0}\}
  \le 2^{n/(18r)} p^{n-n/(18r)}.
\]
Taking into account the number of at most $p^r$ possibilities for
$\mathbf{a}$, the number of singular common zeros $\mathbf{x}$ of
$C_1^*, \ldots, C_r^*$ can therefore be bounded above by
\begin{equation}
\label{nil}
  p^r 2^{n/(18r)} p^{n-n/(18r)}.
\end{equation}
Our aim is to find at least one non-singular zero. Comparing
\eqref{auf_dem} and \eqref{nil} we are led to the inequality
\[
  p^{n-3r}>
  p^r 2^{n/(18r)} p^{n-n/(18r)},
\]
which is satisfied for
\[
  n>72r^2 \frac{1}{1-\log 2/\log p}.
\]
This finishes the proof.

\bibliography{syscubp2}
\bibliographystyle{plain}
\end{document}